\documentclass[12pt]{article}
\usepackage{theorem}
\usepackage{graphics}
\newtheorem{lemma}{Lemma}
\newtheorem{proposition}[lemma]{Proposition}
\newtheorem{theorem}[lemma]{Theorem}
\newtheorem{corollary}[lemma]{Corollary}
{\theorembodyfont{\upshape}}
{\theorembodyfont{\upshape}}
{\theorembodyfont{\upshape}}
{\theorembodyfont{\upshape}}
{\theorembodyfont{\upshape}}
{\theorembodyfont{\upshape}}
{\theorembodyfont{\upshape}}

\newcommand{\R}{{\bf R}}
\newcommand{\C}{{\bf C}}
\newcommand{\rme}{{\rm e}}
\newcommand{\rmd}{\,{\rm d}}

\newcommand{\cB}{{\cal B}}
\newcommand{\cC}{{\cal C}}

\newcommand{\cE}{{\cal E}}
\newcommand{\cF}{{\cal F}}

\newcommand{\cH}{{\cal H}}

\newcommand{\cL}{{\cal L}}

\newcommand{\cT}{{\cal T}}
\newcommand{\sig}{\sigma}
\newcommand{\alp}{\alpha}
\newcommand{\bet}{\beta}
\newcommand{\gam}{\gamma}
\newcommand{\lam}{\lambda}
\newcommand{\Lam}{\Lambda}
\newcommand{\del}{\delta}
\newcommand{\eps}{\varepsilon}
\newcommand{\Gam}{\Gamma}
\newcommand{\kap}{\kappa}
\newcommand{\ome}{\omega}
\newcommand{\Ome}{\Omega}

\newcommand{\Dom}{{\rm Dom}}

\newcommand{\Spec}{{\rm Spec}}

\newcommand{\Ker}{{\rm Ker}}
\newcommand{\Ran}{{\rm Ran}}
\newcommand{\norm}{\Vert}
\renewcommand{\Re}{{\rm Re}\,}
\renewcommand{\Im}{{\rm Im}\,}
\newcommand{\supp}{{\rm supp}}
\newcommand{\Proof}{\underbar{Proof}{\hskip 0.1in}}
\newcommand{\Note}{\underbar{Note}{\hskip 0.1in}}
\newcommand{\hash}{\#}

\newcommand{\Schrodinger}{Schr\"odinger }

\newcommand{\pr}{\prime}
\newcommand{\emp}[1]{{\it #1}}

\newcommand{\nopic}[1]{}

\newcommand{\ops}{one-parameter semigroup}

\newcommand{\be}{\begin{equation}}
\newcommand{\ee}{\end{equation}}

\usepackage{hyperref}
\title{Semi-classical Analysis and Pseudospectra\\
sap2.tex}
\author{E.B. Davies}
\date{10 February 2004}

\begin{document}
\maketitle

\begin{abstract}
We prove an approximate spectral theorem for non-self-adjoint
operators and investigate its applications to second order
differential operators in the semi-classical limit. This leads to
the construction of a twisted FBI transform. We also investigate
the connections between pseudospectra and boundary conditions in
the semi-classical limit.
\end{abstract}

\section{Introduction}

In the last ten years the theory of pseudospectra has developed
rapidly, and has been shown to give substantial insights into the
properties of non-self-adjoint (NSA) matrices and operators,
\cite{boe,EBD3,EBDwebsite,ETwebsite,tre1,tre2}. In this paper we
focus on its applications to second order differential operators.
This involves giving a new and more general definition of
pseudospectra. Our first reason for extending the concept is that
the standard definition does not provide any link with the
geometry of phase space, which is of great importance in the
theory of differential and pseudodifferential operators. By
incorporating the connection into the definitions, we increase the
conceptual clarity and facilitate the analysis of pseudospectra in
those situations in which the semi-classical approximation is
relevant.

The second reason for concentrating on pseudo-eigenfunctions
rather than pseudospectra is that the former are used in
\cite{EBD2} to provide a new method of solving evolution equations
approximately. In several dimensions one could not hope to obtain
sufficient pseudo-eigenfunctions by choosing just one for each
point of the complex plane. Questions of spectral multiplicity
arise just as they do for ordinary spectral theory, and indicate
that a better parametrization is by points in the classical phase
space, not by complex numbers. We plan to use the results of this
paper to extend those of \cite{EBD2} to more general operators.

The paper has three parts. In the first we prove an abstract
approximate spectral theorem for NSA operators. We find a
connection between this and quantization. The second part relates
these ideas to the semi-classical analysis of differential
operators via the semi-classical principal symbol of the operator
and what we call interior pseudo-eigenvectors. Finally we
introduce the concept of boundary pseudo-eigenvectors and describe
how to construct them. We mention that \cite{CT} contains results
relating the boundary and interior pseudospectra of twisted
Toeplitz operators which are parallel to the ones which we obtain
for differential operators. See \cite{DT} for related work on the
wave equation.

\section{An Approximate Spectral Theorem}

In \cite{EBD2} we have shown how to `diagonalize' highly
non-normal operators by using pseudospectra. The diagonalization
is only approximate, but, in spite of this, it may be used to
solve evolution equations efficiently for some quite singular
infinitesimal generators.

In this paper we formulate the underlying theorem at a general
level, in order to make it accessible to a wider audience. All of
the assumptions here are satisfied in the numerical examples
discussed in \cite{EBD2}, as we indicate in the next section. The
ingredients are simple. We suppose that $A$ is a bounded or
closed, unbounded linear operator acting in a separable Hilbert
space $\cH$. We also suppose that $\Lam$ is a multiplication
operator acting in the space $L^p(\Ome, \rmd \ome)$ where $1\leq p
<\infty$; for numerical calculations the simplest choice is $p=2$,
but $p=1$ is more natural for some other purposes. We assume
explicitly that
\[
(\Lam \psi)(\ome)=\sig(\ome)\psi(\ome)
\]
for all $\psi$ in the maximal subdomain of $L^p(\Ome)$, where the
`symbol' $\sig:\Ome \to \C$ of the operator $A$ is a measurable
function and $\rmd \ome$ is a $\sig$-finite measure on $\Ome$. It
is known that the spectrum of the operator $\Lam$ equals the
essential range of $\sig$. We also assume that $E: L^p(\Ome)\to
\cH$ is a bounded linear operator such that
$E(\Dom(\Lam))\subseteq \Dom(A)$ and that
\begin{equation}
\norm AE-E\Lam  \norm <\eps  \label{fundass}
\end{equation}
for a (preassigned, small) $\eps >0$, in the sense that
\begin{equation}
\norm AE\phi-E\Lam\phi  \norm <\eps\norm \phi\norm\label{fundass2}
\end{equation} for all $\phi\in\Dom(\Lam)$.

\begin{theorem} \label{main}
Let $A$ be the generator of a \ops\ $T_t$ acting on $\cH$ and
satisfying
\begin{equation}
\norm T_t\norm\leq M\rme^{\gam t} \label{semibound}
\end{equation}
for all $t\geq 0$. Suppose also that
\[
\Re(\sig(\ome))\leq \gam
\]
for all $\ome\in\Ome$. Then (\ref{fundass}) implies
\begin{equation}
\norm T_tE -E\rme^{\Lam t} \norm \leq    %
\eps tM \rme^{\gam t}\label{estimate}
\end{equation}
for all $t\geq 0$.
\end{theorem}

\Proof Since the operators in (\ref{estimate}) are all bounded it
is sufficient to prove the estimate for all $\phi\in\Dom(\Lam)$.
We then have
\begin{eqnarray*}
\norm T_tE\phi-E\rme^{\Lam t}\phi\norm &= &   \norm \int_0^t
\frac{\rmd}{\rmd s}( T_{t-s}E\rme^{\Lam  s}
 \phi)\, \rmd s \norm\\
&\leq &  \int_0^t \norm  T_{t-s}(AE-E\Lam)\rme^{\Lam  s}
 \phi\norm\, \rmd s \\
&\leq & \int_0^t \norm  T_{t-s}\norm  \eps\norm
 \rme^{\Lam s}\phi\norm\, \rmd s \\
 &\leq & \int_0^t M\rme^{\gam(t-s)} \eps \rme^{\gam s}\norm
\phi\norm\, \rmd s \\
 &= & \eps\norm \phi\norm tM \rme^{\gam t}.
\end{eqnarray*}

If $A$ is a bounded normal operator then the spectral theorem
states that one can find such a representation in which $E$ is
unitary, $\eps =0$ and the essential range of $\Lam$ equals the
spectrum of $A$. The point of Theorem~\ref{main} is that it may be
applied to operators which are far from unitary and in situations
in which the essential range of $\Lam$ is very different from the
spectrum of $A$. The explanation of this relates to pseudospectral
theory.

One might try to develop an `approximate functional calculus'
based upon the above theorem. For example if $T_t=\rme^{At}$ is a
contraction semigroup then under suitable conditions one can prove
an analogue of Theorem~\ref{main} for
$T_{\alp,t}=\rme^{-(-A)^{\alp}t}$ when $0<\alp<1$; see
\cite{EBD2}.

In order to compare Theorem~\ref{main} with the results in
\cite{EBD2} one needs to approximate $E$ by an operator $E^\pr$
whose range is not contained in $\Dom(A)$.

\begin{corollary} \label{Eprime}
If in addition to the previous assumptions one has
$\norm E-E^\pr\norm<\eps$ then
\begin{equation}
\norm T_tE^\pr\phi -E^\pr\rme^{\Lam t} \phi\norm \leq    %
 \eps\norm \phi\norm (1+M+ tM) \rme^{\gam t}\label{estimate2}
\end{equation}
for all $\phi\in L^p(\Ome)$ and all $t\geq 0$.
\end{corollary}

\Proof This follows directly from
\begin{eqnarray*}
\norm T_tE^\pr\phi -E^\pr\rme^{\Lam t} \phi\norm &\leq & \norm
T_tE\phi
-E\rme^{\Lam t} \phi\norm\\
&&+ \norm T_t (E-E^\pr)\phi\norm+ \norm (E-E^\pr)\rme^{\Lam
t}\phi\norm.
\end{eqnarray*}

The following modification of Theorem~\ref{main} assumes that one
is given $f\in \cH$ and wishes to approximate $T_tf$.

\begin{corollary}
If $f\in\cH$ then under the conditions of Theorem~\ref{main}
\[
\norm T_t f-E\rme^{\Lam t}\phi\norm\leq \norm f-E\phi\norm
M\rme^{\gam t} +\eps\norm \phi \norm tM
 \rme^{\gam t}
\]
for all $\phi\in L^p(\Ome)$ and $t\geq 0$.
\end{corollary}

\Proof We have
\[
\norm T_t f-E\rme^{\Lam t}\phi\norm\leq \norm T_t (f-E\phi)\norm+
\norm T_t E\phi-E\rme^{\Lam t}\phi\norm
\]
each of which is straightforward to estimate.

The above results can only be useful if $M$, $t$ and $\gam$ are of
order $1$. There also has to exist $\phi$ such that $\norm
f-E\phi\norm$ and $\eps\norm \phi\norm$ are both small. One cannot
simply put $\phi=E^{-1}f$, since $E$ need not be surjective or
invertible.

If $p=2$, the standard way of solving this problem is to minimize
the functional
\begin{equation}
\cE(\phi)=\norm f-E\phi\norm^2+\del \norm \phi\norm^2  \label{fnl}
\end{equation}
for a suitable value of the regularization parameter $\del >0$;
see \cite{han}. This is achieved in the numerical context by
putting
\[
\phi=\tilde{E}\backslash (f\oplus 0)
\]
where $\tilde{E}:L^2(\Ome)\to \cH\oplus L^2(\Ome)$ is defined by
\begin{equation}
\tilde{E}\phi=E\phi\oplus \del^{1/2}\phi.\label{Etilde}
\end{equation}
We include the proof of the following well-known proposition for
completeness.

\begin{proposition}\label{fnltheorem}
If $p=2$, the minimum of (\ref{fnl}) is achieved for $\phi=F_\del
f$, where
\begin{equation}
F_\del=(E^\ast E+\del I)^{-1}E^\ast\label{Fdeldef}
\end{equation}
satisfies $\norm F_\del\norm\leq\del^{-1/2}$. Moreover $\norm
EF_\del\norm\leq 1$ for all $\del>0$. One has
\begin{equation}
\lim_{\del\to 0}EF_\del f=f \label{limass}
\end{equation}
for all $f\in\cH$ if and only if $\Ran(E)$ is dense in $\cH$.
\end{proposition}

\Proof The first statement depends upon a routine variational
calculation. For the second we observe that
\[
\norm (E^\ast E+\del I)^{-1}E^\norm\leq ab
\]
where
\[
a=\norm (E^\ast E+\del I)^{-1/2}\norm \leq \del^{-1/2}
\]
and
\begin{eqnarray*}
b^2&=&\norm (E^\ast E+\del I)^{-1/2}E^\ast\norm^2\\
 &=&\norm
(E^\ast E+\del I)^{-1/2}E^\ast .\, E(E^\ast E+\del
I)^{-1/2}\norm\\
&\leq & \norm (E^\ast E+\del I)^{-1/2}(E^\ast E+\del I)(E^\ast
E+\del I)^{-1/2}\norm\\
&=& 1.
\end{eqnarray*}
This calculation also implies that
\begin{eqnarray*}
\norm EF_\del\norm&=&\norm E(E^\ast E+\del
I)^{-1/2}.\, (E^\ast E+\del I)^{-1/2}E^\ast \norm\\
&=&\norm (E^\ast E+\del I)^{-1/2}E^\ast\norm^2\\
&\leq & 1.
\end{eqnarray*}
Since $\Ran(EF_\del)\subseteq \Ran(E)$, (\ref{limass}) implies
that $\Ran(E)$ is dense. If $\Ran(E)$ is dense then the uniform
boundedness just proved implies that (\ref{limass}) holds for all
$f\in\cH$ if it holds whenever $f=E\phi$ for some $\phi\in
L^2(\Ome)$. In this case let $P$ denote the orthogonal projection
onto the closure of the range of $E^\ast E$. Since
$\Ker(E)=\Ker(E^\ast E)$, we may assume without loss of generality
that $P\phi=\phi$. We have
\begin{eqnarray*}
\lim_{\del \to 0}EF_\del f&=&\lim_{\del \to 0}E(E^\ast E+\del
I)^{-1}E^\ast E\phi\\
&=& EP\phi=E\phi=f
\end{eqnarray*}
by applying the spectral theorem to the non-negative self-adjoint
operator $E^\ast E$.

Using Proposition~\ref{fnltheorem} one may ensure that $\eps\norm
\phi\norm$ is small by choosing $\del$ appropriately. Even if $E$
has dense range, one cannot ensure that $\norm f-E\phi\norm$ is
small for some particular $\del >0$ without further conditions.
One has either to make the a priori assumption that $f$ lies in
some subspace of well-approximable vectors, or observe a
posteriori for particular choices of $f$ and $\del$ that the
minimizing $\phi$ does indeed make this quantity small enough for
the application intended.


\section{The Connection with Pseudospectra}

Given $\eps >0$, the $\eps$-pseudospectrum of the closed operator
$A$ is defined by
\[
\Spec_\eps(A)=\{ z: \norm Af-zf\norm <\eps \norm f\norm \mbox{ for
some } f\in \Dom(A)\}.
\]
Pseudospectral ideas lie at the core of this paper, and we refer
to \cite{boe,EBD3,EBDwebsite,ETwebsite,tre1,tre2} for background
material on this subject. The following theorem is valid for all
$p\in[1,\infty)$, but its main application is for $p=1$. Indeed we
conjecture that if $p=2$ the first condition on $E$ can only hold
if $E$ is isometric. In the following theorem $P_U$ denotes the
operator of multiplication by the characteristic function of the
set $U$, always assumed to be measurable.

\begin{theorem} \label{pspectral}
Suppose that $1\leq p <\infty$, $\norm EP_U\norm =1$ for all
subsets $U$ of $\Ome$ with positive measure, and $\norm
AE-E\Lam\norm<\eps$. Then
\[
\Spec(\Lam)\subseteq \Spec_\eps(A).
\]
\end{theorem}

\Proof Let $\bet\in\Spec(\Lam)$. We choose $\del>0$ such that
\[
\eps^\pr:=\norm AE-E\Lam\norm+\del <\eps
\]
and put
\[
U=\{ \ome\in\Ome:|\lam(\ome)-\bet|<\del\}.
\]
If $\phi$ has support in $U$ then
\begin{eqnarray*}
\norm AE\phi-\bet E\phi\norm &\leq& \norm (AE-E\Lam)\phi \norm
+\norm E(\Lam\phi-\bet\phi)\norm\\
&\leq & \eps^\pr\norm \phi\norm.
\end{eqnarray*}
Therefore
\begin{eqnarray*}
\inf\{\norm Af-\bet f\norm/\norm f\norm:0\not= f\in \cH\}
&\leq&\inf\{\norm AE\phi-\bet E\phi\norm/\norm
E\phi\norm:0\not=\phi\in L^p(U)\}\\ &\leq& \eps^\pr\inf\{\norm
\phi\norm/\norm
E\phi\norm:0\not=\phi\in L^p(U)\}\\
&=&\eps^\pr <\eps.
\end{eqnarray*}
This implies that $\bet\in\Spec_\eps(A)$.

\begin{theorem} \label{mainseq}
If $p=1$ then the conditions of Theorem~\ref{pspectral} are
equivalent to the following statements. There exists $\eps^\pr >0$
and a set $N$ of zero measure, such that for each
$\ome\in\Ome\backslash N$ there is a unit vector
$e_\ome\in\Dom(A)$ which depends measurably on $\ome$ and
satisfies
\[
\norm Ae_\ome-\sig_\ome e_\ome\norm\leq \eps^\pr<\eps
\]
for all $\ome\in\Ome\backslash N$, where $\sig(\ome)\in\C$.
\end{theorem}

\Proof The passage from the assumptions of Theorem~\ref{pspectral}
to the statements of this theorem is justified by using
\cite[Theorem VI.8.6]{DS}.

If $\phi$ lies in the maximal domain of $\Lam$ then under the
assumptions of Theorem~\ref{mainseq}
\begin{equation}
\norm AE\phi- E\Lam\phi\norm\leq \int_\Ome |\phi(\ome)|\,\norm
Ae_\ome-\sig_\ome e_\ome\norm\,\rmd\ome \leq
\eps^\pr\norm\phi\norm.\label{comestim}
\end{equation}
Hence $\norm AE-E\Lam\norm \leq \eps^\pr <\eps$. The calculations
involved would be easy to justify if one only had to deal with
finite sums, or if $A$ and $\Lam$ were bounded, but in general
they use limiting processes to define the integrals. Commuting $A$
and $\Lam$ with these limiting processes is justified by the
following lemma.

\begin{lemma}
Let $A$ be a closed linear operator with domain in a Banach space
$\cB$ and range in a Hilbert space $\cH$. Let $c>0$,
$f_n\in\Dom(A)$, $\norm f_n-f\norm\to 0$, $\norm g_n-g\norm\to 0$
and $\norm Af_n-g_n\norm\leq c$ for all $n$, then $f\in\Dom(A)$
and $\norm Af-g\norm\leq c$.
\end{lemma}

If $\Ome$ has finite measure $|\Ome|$, then $L^2(\Ome)$ is
continuously embedded in $ L^1(\Ome)$, and all of the theorems of
Section 1 hold under the present hypotheses. In the numerical
applications of \cite{EBD2} the space $\Ome$ is taken to be the
finite set $\{1,.., N\}$ and $\rmd \ome$ is the counting measure.
Given  unit pseudo-eigenvectors $e_n\in\cH$ of $A$ for $1\leq
n\leq N$, we have
\begin{equation}
E\phi=\sum_{n=1}^N \phi_ne_n.\label{Ephi}
\end{equation}
There is no requirement that the vectors should be linearly
independent, and indeed in some of the examples studied in
\cite{EBD2} they are taken from an overcomplete infinite sequence
$\{e_n\}_{n=1}^\infty$. Equivalently the operator $E$ need not be
invertible, or may have a large condition number.


\section{Quantization}

In this section we make some general comments about the
relationship between our previous results and the notion of
quantization.

Let $\Ome$ be a second countable locally compact Hausdorff space,
and let $\rmd\ome$ be a regular Borel measure on $\Ome$ with
support equal to $\Ome$. Let $\cH$ be a separable Hilbert space
and let $e:\Ome\to \cH$ be a continuous function. We define
$E:C_c(\Ome)\to\cH$ by
\[
E\phi=\int_\Ome \phi(\ome) e_\ome\,\rmd \ome.
\]
The following are well-known and elementary.

\begin{lemma} The operator $E$ extends to a bounded linear
operator $E_1:L^1(\Ome,\rmd\ome)\to \cH$ if and only if
$\ome\to\norm e_\ome\norm$ is a bounded function, in which case
\[
\norm E_1\norm=\sup\{\norm e_\ome\norm:\ome\in\Ome\}.
\]
The operator $E$ extends to a bounded linear operator
$E_2:L^2(\Ome,\rmd\ome)\to \cH$ if and only if
\begin{equation}
 \int_\Ome |\langle f,e_\ome\rangle|^2\,\rmd \ome \leq
c^2\norm f\norm^2\label{L2bound}
\end{equation}
for some $c\geq 0$ and all $f\in\cH$, in which case $\norm
E_2\norm$ is the smallest such constant $c$. The operator
$E^\ast:\cH\to C(\Ome)$ is an isometry from $\cH$ into
$L^2(\Ome,\rmd\ome)$ if and only if
\begin{equation}
\int_\Ome |\langle f,e_\ome\rangle|^2\,\rmd \ome =\norm
f\norm^2\label{overcomplete}
\end{equation}
for all $f\in\cH$.
\end{lemma}

Families of vectors $\{e_\ome\}_{\ome\in\Ome}$ satisfying
(\ref{overcomplete}) are also called continuous resultions of the
identity and have played an important part in group representation
theory and quantum mechanics for many decades. For their
connection with coherent state theory and the Bargman transform
see \cite[ch. 8]{QTOS} and \cite[ch. 3]{mar}. If
(\ref{overcomplete}) holds then $E\backslash f=E^\ast f$ for all
$f\in\cH$. but this is not the case under the assumption
(\ref{L2bound}), which is more relevant to this paper.

Given a function $f\in C_c(\Ome)$ we define the multiplication
operator $M_f$ by $M_f\phi=f\phi$ where $\phi\in L^p(\Ome)$ for
some $p$. We define the quantization of the function $f$ to be the
operator $Q(f)=EM_f E^\ast$ on $\cH$. We may also write
\[
Q(f)=\int_\Ome f(\ome) P_{e_\ome}\,\rmd\ome
\]
where $P_a \psi=\langle \psi,a\rangle a\,$; see, for example,
\cite[sect. 8.5]{QTOS}. The following lemma is also standard.

\begin{lemma}
If $f\geq 0$ then $Q(f)\geq 0$. If $E$ is bounded from $L^1(\Ome)$
to $\cH$ then $Q$ extends to a bounded linear operator from
$L^1(\Ome)$ to the space $\cT(\cH)$ of trace class operators on
$\cH$. If $E$ is bounded from $L^2(\Ome)$ to $\cH$ then $Q$
extends to a bounded linear operator from $L^\infty(\Ome)$ to the
space $\cL(\cH)$ of bounded operators on $\cH$. Given
(\ref{overcomplete}), or equivalently $EE^\ast=1$, we have
$Q(1)=1$.
\end{lemma}

In quantum theory it is commonplace to refer not to the operator
$Q$ but to the positive-operator-valued measure $A(U):=EM_{\chi_U}
E^\ast$ where $\chi_U$ is the characteristic function of the
measurable set $U$ of $\Ome$. The formula
\[
Q(f)=\int_\Ome f(\ome)\,A(\rmd \ome)
\]
implements a one-one correspondence between the two definitions;
see \cite[Lemma 3.1.2]{QTOS}. If $EE^\ast=1$ then $A(\Ome)=1$ and
$A(\cdot)$ is called a generalized observable; for a systematic
study of POV measures and their relation to coherent states see
\cite[Ch. 3]{QTOS} or \cite{hol}. See \cite{HS} for more recent
references and a connection with subnormal operators.

The difference between this method of quantization and the
approach of this paper is now clear. Instead of studying
$Q(f)=EM_f E^\ast$, we would like to study $S(f)=EM_f E^{-1}$. If
this were possible $f\to S(f)$ would be an algebra homomorphism
from $L^\infty(\Ome)$ to $\cL(\cH)$.  Since $E$ is not invertible
in general we compromise by studying $EM_f F_\del$, where the
regularized inverse $F_\del$ is given by (\ref{Fdeldef}) and
$\del>0$ is chosen small enough to yield numerically valuable
results but not so small that the computational algorithms become
unreliable.

The operator $E$ which we have considered above has much in common
with the Fourier-Bros-Iagolnitzer (FBI) transform as defined in
\cite[ch. 3]{mar}. See also \cite[Ch.3]{QTOS}, where the
connection with the Wigner distribution and applications to
quantum theory are explained. In Section \ref{transform} we define
a distorted FBI transform; the distortions are introduced to adapt
the transform to a given differential operator, and involve
replacing the Gaussian states used in the definition of the FBI
transform by pseudo-eigenfunctions of the operator.

\section{The Connection with Semi-classical Analysis}

Before describing the connection of the above ideas with
semi-classical analysis, we generalize the notion of
pseudospectra. Following \cite{EBDwebsite,HT,lav,TH}, we define
the (generalized) pseudospectra of a family of closed operators
$\{A_\ome\}_{\ome\in\Ome}$ acting from dense domains
$\Dom(A_\ome)$ in a Banach space $\cB$ to another Banach space
$\cC$ to be the sets
\[
\Spec_\eps(A)=\{ \ome: \norm A_\ome f\norm<\eps \norm f\norm
\mbox{ for some $f\in\Dom(A_\ome)$}\}
\]
where $\eps >0$. We have
\[
\Spec_\eps(A)\cup \Spec_\eps(A^\ast)=\Spec(A)\cup\{ \ome:\norm
A_\ome^{-1}\norm >\eps^{-1}\}
\]
where $\Spec(A)$ is defined to be the set of $\ome$ for which
$A_\ome$ is not invertible. If $\dim(\cB)=\dim(\cC)<\infty$ then
\[
\Spec(A)\subseteq\Spec_\eps(A)=\Spec_\eps(A^\ast)
\]
for all $\eps>0$. If $\dim(\cB) <\dim(\cC)<\infty$ then
\[
\Spec_\eps(A^\ast)=\Spec(A)=\Ome
\]
for all $\eps >0$. The proof of the following lemma may be found
in \cite{HT}.

\begin{lemma} One has $\ome\in\Spec_\eps(A)$ if and only if there
exists a bounded operator $D:\cB\to\cC$ such that $\norm D \norm
<\eps$ and
\[
\Ker(A(\ome)+D)\not= \{0\}.
\]
\end{lemma}

Given a differential or pseudodifferential operator $L_h$ with
domain $C^\infty_c(X)$, where $X$ is a region in $\R^N$ and $h>0$,
we define the operator family
\[
A_{h,u,\xi}:C^\infty_c(X)\subseteq L^2(X) \to L^2(X,\C^{2N+1}),
\]
by
\begin{equation}
A_{h,u,\xi} f=(Q^jf-u^jf,P_jf-\xi_jf,L_hf-\sig(u,\xi)f)\label{Ah}
\end{equation}
where $(Q^jf)(x)=x^jf(x)$ and $(P_jf)(x)=-ih\partial_jf(x)$. In
these equations we assume that $u\in X$, $\xi\in\R^N$, $1\leq
j\leq N$ and $\sig(u,\xi)$ is the semiclassical principal symbol
of the operator $L_h$, as defined below. It follows directly from
the definitions that $\norm A_{h,u,\xi}f \norm <\eps \norm f\norm$
implies
\begin{eqnarray*}
\norm Q^jf -u^jf\norm&< &\eps\,\norm
f\norm\\
\norm P_jf -\xi_jf\norm&<&\eps\,\norm
f\norm\\
\norm L_hf -\sig(x,\xi)f\norm&< &\eps\,\norm f\norm
\end{eqnarray*}
where $1\leq j \leq N$. It is known that the pseudospectra
converge to fill a certain set $\sig(\Lam)$ if $h\to 0$ and
$\eps\to 0$ simultaneously at suitable rates; see Section 7 for
details. Even in one space dimension a point in $\sig(\Lam)$ may
be the image of more than one point in $\Lam$, so $\sig(\Lam)$ may
have hidden structure as a subset of $\C$. This observation
applies with less precision to the numerically determined
pseudospectra for fixed $h> 0$ and $\eps
>0$.

The extension of the above ideas to a manifold $X$ needs some
care, since the full symbol $\sig_h(u,\xi)$ is not an invariant
object in general. It is shown in \cite{saf} that one can resolve
these problems if the manifold is provided with a linear
connection, as happens if it is Riemannian. The symbol
$\sig_h(u,\xi)$ is then definable as a function on the cotangent
bundle $T^\ast X$ and $\Lam$ is a certain subset of $T^\ast X$. We
do not actually need the full symbol for our problem: its
semiclassical limit is sufficient. The semiclassical principal
symbol is given by
\[
\sig(u,\xi)=\lim_{h\to 0} \sig_h(u,h^{-1}\xi)
\]
and is an invariant quantity, i.e. as a function on the cotangent
bundle $T^\ast X$ it does not depend on the choice of local
coordinates.

The following alternative definition of the semiclassical
principal symbol of $L_h$ makes its invariant character clear.
Suppose that $u\in X$ and $\xi$ is a cotangent vector at $u$. Let
$f$ be any smooth function on $X$ such that $\rmd f(u)=\xi$. Then
\[\sig(u,\xi)=\left\{\lim_{h\to 0}
\rme^{-ih^{-1}f}L_h\left(\rme^{ih^{-1}f}\right)\right\}(u).
\]


\section{The Semiclassical Spectrum}

The theory which we shall describe can be developed at several
levels of generality, and in this section we consider only second
order differential operators acting on $\R^N$.

Given $h>0$, let $L_h$ denote the operator
\[
(L_hf)(x)=-h^2a^{j,k}_h(x)\partial_{j,k}f(x)-ihb^j_h(x)\partial_jf(x)+c_h(x)f(x)
\]
acting on functions $f:\R^N\to \C$, where $a,b,c$ are sufficiently
regular functions whose values are respectively matrices, vectors
and scalars with complex-valued entries, and we use the standard
summation convention. Under conditions which we shall impose the
domain of $L_h$ will contain $C^\infty_c(\R^N)$. All
considerations in this paper are local, so no growth bounds at
infinity on the coefficients are needed. We allow the coefficients
to be $h$-dependent so that the class of differential operators is
invariant under local changes of coordinates. The semiclassical
principal symbol of this operator is the complex-valued function
\begin{equation}
\sig(u,\xi)=
a^{j,k}_0(u)\xi_j\xi_k+b^j_0(u)\xi_j+c_0(u)\label{scps}
\end{equation}
in which we take $u,\xi$ to be \emp{real} vectors in $\R^N$.

Given $(u,\xi)\in \R^N\times \R^N$ we are interested in finding
localized approximate eigenfunctions for the operator $L_h$. We
require that they become asymptotically exact as $h\to 0$.

Our first theorem provides the motivation for defining the
semi-classical spectrum of $L_h$ to be the set
$\sig(\R^N\times\R^N)$.

\begin{theorem} \label{two} Suppose that $a^{j,k}_h(x)$,
$b^j_h(x)$ and $c_h(x)$ are all locally Lipschitz continuous in
both $x\in\R^N$ and $h\in [0,1]$. Then for every $u\in\R^N$,
$\xi\in\R^N$ and $h\in(0,1]$ there exists $f_h\in
C^\infty_c(\R^N)$ such that
\begin{eqnarray}
\norm f_h\norm_2&=& c>0\label{spec1}\\
\norm Q^jf_h-u^jf_h\norm_2&=& O(h^{1/2})\label{spec2}\\
\norm P_jf_h-\xi_jf_h\norm_2&=& O(h^{1/2})\label{spec3}\\
\norm L_hf_h-\sig(x,\xi)f_h\norm_2&=& O(h^{1/2})\label{spec4}
\end{eqnarray}
as $h\to 0$, for all $1\leq j\leq N$.
\end{theorem}

\Proof Let $\phi$ be a non-negative $C^\infty$ function on $\R^N$
which equals $1$ if $|x|\leq 1$ and $0$ if $|x|\geq 2$. Given
$(u,\xi)\in\R^N\times \R^N$, $h >0$ and $\alp=1/2$  define
\[
f_h(x)=h^{-N\alp/2}\rme^{ih^{-1}\xi\cdot x}\phi(h^{-\alp}(x-u)).
\]
The first three statements of the theorem are routine
verifications performed by the same method as follows.

We verify (\ref{spec4}) by using the expansion
\[
L_hf_h-\sig(u,\xi)f_h=g+r_1+r_2+r_3+r_4
\]
where
\begin{eqnarray*}
g(x)&=& \{ a^{j,k}_h(x)-a^{j,k}_0(u)\} \xi_j\xi_k f_h(x)\\
&& +\{ b^j_h(x)-b^j_0(u)\}\xi_jf_h(x)+\{c_h(x)-c_0(u)\}f_h(x)
\end{eqnarray*}
and
\begin{eqnarray*}
r_1&=&-ih^{1-\alp-N\alp/2}a^{j,k}(x)\xi_j\rme^{ih^{-1}\xi\cdot
x}\phi_k(h^{-\alp}(x-u))\\
r_2&=&-ih^{1-\alp-N\alp/2}a^{j,k}(x)\xi_k\rme^{ih^{-1}\xi\cdot
x}\phi_j(h^{-\alp}(x-u))\\
r_3&=&h^{2-2\alp-N\alp/2}a^{j,k}(x)\rme^{ih^{-1}\xi\cdot
x}\phi_{j,k}(h^{-\alp}(x-u))\\
r_4&=&-ih^{1-\alp-N\alp/2}b^{j}(x)\xi_j\rme^{ih^{-1}\xi\cdot
x}\phi_j(h^{-\alp}(x-u)).
\end{eqnarray*}
In these identities the subscripts on $\phi$ denote partial
derivatives. The Lipschitz assumptions on the coefficients of
$L_h$ and the fact that the support of $f_h$ has diameter of order
$h^\alp$ imply that
\[
\norm g\norm_2 =O(h^\alp)
\]
as $h\to 0$. We also have $\norm r_j\norm_2 =O(h^{1-\alp})$ for
$j=1,2,4$ and $\norm r_3\norm_2 =O(h^{2-2\alp})$. The overall
error is minimized by putting $\alp=1/2$.


\section{Constructing the Interior Pseudospectra}

The material in this section is based upon the fact that if the
coefficients are sufficiently smooth then the estimate
(\ref{spec4}) can be greatly improved by a suitable choice of
$f_h$. In the language of Section 5 we replace (\ref{Ah}) by
\begin{equation}
A_{h,u,\xi} f=(Q^jf-u^jf,P_jf-\xi_jf,h^{-n}\{
L_hf-\sig_0(u,\xi)f\} )\label{Ah2}
\end{equation}
where $n>0$. The size of $n$ depends upon the smoothness of the
coefficients, which for simplicity we assume to be $C^\infty$. The
pseudospectral estimate $\norm A_{h,u,\xi}f\norm <\eps\norm f
\norm$ then implies
\begin{eqnarray*}
\norm Q^jf -u^jf\norm&< &\eps\,\norm
f\norm\\
\norm P_jf -\xi_jf\norm&<&\eps\,\norm
f\norm\\
\norm L_hf -\sig(x,\xi)f\norm&< &h^{n}\eps\,\norm f\norm
\end{eqnarray*}
where $1\leq j \leq N$. We repeat the calculations of
\cite{EBD5,EBD6} for a more general second order ordinary
differential operator for completeness. The extension to
pseudo-differential operators in higher dimensions,
\cite{dsz,zworski}, cannot be formulated in exactly the same
manner: there can be infinitely many different
pseudo-eigenfunctions associated with a point in phase space, and
the correct parametrization of these is not obvious. We assume
that
\[
(L_hf)(x)=-h^2a(x)f^{\pr\pr}(x)-ihb(x)f^\pr(x)+c(x)f(x)
\]
so that the semiclassical principal symbol is
\[
\sig(u,\xi)=a(u)\xi^2+b(u)\xi +c(u).
\]
We assume ellipticity, in other words that $a(x)\not= 0$ for all
$x\in\R$. Given $u,\xi\in\R$, we put
\begin{equation}
f(u+s)=h^{-1/4}\chi(s)\exp(\psi(s))\label{fdef}
\end{equation}
for all $s\in\R$, where $\chi\in C_c^\infty$ satisfies $\chi(s)=1$
if $|s|\leq\del/2$ and $\chi(s)=0$ if $|s|\geq \del$, and $\del
>0$ must be small enough; see the proof of Lemma~\ref{asympt}. We
assume that
\begin{equation}
\psi(s)=\sum_{m=-1}^n h^m\psi_m(s)\label{psidef}
\end{equation}
for some integer $n\geq -1$. This is a non-standard form of the
JWKB expansion, and has the feature that the function $f$ does not
vanish within the interval of interest. A direct computation shows
that
\begin{equation}
L_h f-\sig(u,\xi)f= \left(\sum_{m=0}^{2n+2}h^m\phi_m\right)f
+\rm{Rem}\label{pseudoexpand}
\end{equation}
where ${\rm Rem}=O(h^\infty)$ as $h\to 0$ under the conditions
which we impose below. Also
\begin{eqnarray*}
\phi_0(s)&=&-a(u+s)(\psi_{-1}^\pr(s))^2-ib(u+s)\psi_{-1}^\pr(s)
+c(u+s)\\
&&-a(u)\xi^2-b(u)\xi-c(u).
\end{eqnarray*}
Assuming ellipticity, that is $a(x)\not= 0$ for all $x\in\R$, the
eikonal identity $\phi_0=0$ implies
\[
\psi_{-1}(s)=i\int_{v=0}^s\left\{
-\frac{b(u+v)}{2a(u+v)}+\sqrt{w(u,\xi,v)}\right\}\,\rmd v
\]
where
\[ w(u,\xi,v)= \frac{a(u)\xi^2}{a(u+v)} + \frac{b(u)\xi}{a(u+v)}
+\frac{b(u+v)^2}{4a(u+v)^2} +\frac{c(u)-c(u+v)}{a(u+v)}.
\]
We take the branch of the square root which equals
$\xi+b(u)/2a(u)$ at $v=0$. The condition (\ref{nontriv}) implies
that $\partial \sig /\partial \xi \not= 0$ and hence that
$w(u,\xi,0)$ is non-zero; this implies that $w(u,\xi,v)\not= 0$
for all small enough $v$; and hence that the square root is
uniquely determined for all such $v$ by the requirement of
continuity.

Writing $\psi_{-1}(s)=i\xi s+ks^2/2+O(s^3)$ for some $k\in\C$, we
then obtain
\[
-ik\{2a(u)\xi +b(u)\}+a^\pr(u)\xi^2+b^\pr(u)\xi+c^\pr(u)=0.
\]
The requirement that $\Re(k)<0$ may be rewritten in the form
\[
\Im\left( \frac{\partial\sig}{\partial
u}\frac{\partial\overline{\sig}}{\partial \xi}\right) <0
\]
and then in the form $(u,\xi)\in\Ome$ where
\begin{equation}
\Ome =\{ (u,\xi):\{\sig_1,\sig_2\}>0\}\label{nontriv}
\end{equation}
and
\[
\{\sig_1,\sig_2\}:=\frac{\partial\sig_1}{\partial
u}\frac{\partial\sig_2}{\partial
\xi}-\frac{\partial\sig_1}{\partial
\xi}\frac{\partial\sig_2}{\partial u}.
\]
and $\sig_1=\Re(\sig)$, $\sig_2=\Im(\sig)$. In examples one may
find that $\Ome$ is not connected. If it has components $\Lam_j$
then $\sig(\Ome_j)$ may overlap. The multiplicity of a point $z\in
\sig(\Ome)$ may be defined by
\[
m_L(z)=\hash \{ (u,\xi)\in\Ome:\sig(u,\xi)=z\}.
\]

If the coefficients of $L_h$ are smooth then for any choice of $n$
one may choose $\psi_0,...,\psi_n$ so that
$\phi_1=...=\phi_{n+1}=0$. This is achieved as follows. If $1\leq
m \leq n$ then
\[
\phi_{m+1}=(-2a\psi_{-1}^\pr-ib)\psi_m^\pr+F_m(\psi_{-1},...,\psi_{m-1}).
\]
It follows from (\ref{nontriv}) that $2a\psi_{-1}^\pr+ib\not= 0$
if $s=0$, and hence that it is non-zero for all small enough $s$.
If we define $\psi_m$ by
\[
\psi_m(s)=\int_0^s\frac{F_m(\psi_{-1},...,\psi_{m-1})}{2a\psi_{-1}^\pr+ib}
\,\rmd v
\]
Then $|\psi_m(s)|\leq c_m |s|$ and $\phi_{m+1}(s)=0$ for all small
enough $s$. On making these choices we obtain a
pseudo-eigenfunction $f$, depending on $h,n,u$ and $\xi$, for
which $ L_h f -\sig(u,\xi) f= O(h^{n+2})$ as $h\to 0$.

The proof of Theorem~\ref{mainasympt} below is facilitated by
introducing the scale of spaces $\cE^\gam$, consisting of all
functions which can be written as finite sums of functions of the
form $g(s)=h^{\alp-1/4} s^\bet \rho(s)\exp\{\psi(s)\}$ where
$\psi$ is given by (\ref{psidef}), $\rho \in C^\infty$ has support
in $[-\del,\del]$, $\alp\in\R$, $\bet\in\{0,1,2,...\}$ and
$2\alp+\bet\geq \gam$. Putting
$\cE^\infty=\cap_{\gam\in\R}\cE^\gam$ we see that if, in addition
to the above assumptions, $\rho$ vanishes in some neighbourhood of
$0$, then $g\in \cE^\infty$.

\begin{lemma} \label{asympt}
If $\del >0$ is small enough and $g\in \cE^\gam$ then there exists
$c$ such that
\[
\norm g\norm\leq c h^{\gam/2}
\]
for all $0<h\leq 1$.
\end{lemma}

\Proof It is sufficient to consider the case in which $g$ is one
of the terms of the form assumed in the definition of $\cE^\gam$.
One may rewrite $|h^{1/4 -\alp}g(s)|^2$ in the form
$s^{2\bet}G(s)\exp\{-h^{-1}s^2F(s)\}$ where
$F(s)=-2\Re(\psi_{-1}(s)/s^2)$ is a positive continuous function
on $[-\del,\del]$ if $\del >0$ is small enough and $G$ is a
continuous function on $[-\del,\del]$. By Laplace's method we have
\[
\int_{-\del}^\del s^{2\bet}G(s)\exp\{-h^{-1}s^2F(s)\}\,\rmd s\sim
ch^{({2\bet}+1)/2}
\]
as $h\to 0+$, where
\[
c= \frac{G(0)\Gam(({2\bet}+1)/2)}{ F(0)^{({2\bet}+1)/2}}.
\]
The statement of the lemma follows immediately.

\begin{theorem}\label{mainasympt}
If the coefficients of $L_h$ are $C^\infty$ and $(u,\xi)$ lies in
the set $\Ome$ defined by (\ref{nontriv}), then for every positive
integer $n$ there exist functions $f\in C^\infty_c$ depending on
$h,n,u,\xi$ such that
\begin{eqnarray}
\lim_{h\to 0}\norm f\norm &=& c>0 \label{eq1}\\
\norm Q{f}-u{f}\norm&=& O(h^{1/2})\label{eq2}\\
\norm P{f}-\xi{f}\norm&=& O(h^{1/2})\label{eq3}\\
\norm L_h{f}-\sig(u,\xi){f}\norm&=& O(h^{n+2})\label{eq4}
\end{eqnarray}
as $h\to 0$.
\end{theorem}

\Proof We define $f$ by (\ref{fdef}) and observe that $f\in
\cE^0$. The asymptotic formula (\ref{eq1}) follows by the method
of proof of Lemma~\ref{asympt}. We next observe that $Qf-uf\in
\cE^1$ so (\ref{eq2}) follows from Lemma~\ref{asympt}.

We have
\[
Pf-\xi f = \mu_1+\mu_2+\mu_3
\]
where
\begin{eqnarray*}
\mu_1 &=& -ih^{-1/4}\{\psi_{-1}^\pr(s)-i\xi
\}\chi(s) \exp\{\psi(s)\}\in \cE^1\\
\mu_2 & = &  -ih^{3/4}\{ \sum_{m=0}^n h^m\psi_m^\pr(s) \}\chi(s)\exp\{\psi(s)\}\in \cE^2\\
\mu_3 & = & -ih^{3/4}\chi^\pr(s)\exp\{\psi(s)\}\in \cE^\infty.
\end{eqnarray*}
Therefore $ Pf-\xi f \in \cE^1 $ and (\ref{eq3}) follows using
Lemma~\ref{asympt}.

Since $\phi_m=0$ for $0\leq m\leq n+1$ it follows from
(\ref{pseudoexpand}) that
\begin{eqnarray*}
L_h f-\sig(u,\xi)f&=& \left(\sum_{m=n+2}^{2n+2}h^m\phi_m\right)f
+O(h^\infty)\\
&\in & \cE^{2n+4}
\end{eqnarray*}
This implies (\ref{eq4}) by Lemma~\ref{asympt}.

\Note The orders of magnitude of the errors in (\ref{eq2}) and
(\ref{eq3}) cannot \emp{both} be reduced by a different choice of
the function $f$, because of the uncertainty principle.

The following lemma shows that one can approximate the
pseudo-eigenfunction by a Gaussian expression.

\begin{lemma}\label{gaussapp}
We have
\[
\norm f-g\norm\leq ch^{1/2}
\]
as $h\to 0$, where
\[
g(u+s)=h^{-1/4}\exp\{h^{-1}(i\xi s+ks^2/2)\}.
\]
\end{lemma}

\Proof Since $g-\chi g=O(h^\infty)$ we have to estimate the $L^2$
norm of
\[
h^{-1/4}\chi(s)\left(\exp\{\psi(s)\}-\exp\{h^{-1}(i\xi
s+ks^2/2)\}\right).
\]
By virtue of the bound
\[
|\rme^{-a}-\rme^{-b}|\leq |a-b|\rme^{-\min(\Re(a),\Re(b))}
\]
this is dominated by the absolute value of
\[
\mu(s)=h^{-1/4}\{ \psi(s)- h^{-1}( i\xi
s+ks^2/2)\}\chi(s)\exp\{-h^{-1}cs^2\}
\]
for some $c>0$. In the following calculations we define
$\tilde{\cE}^\gam$ in the same way as $\cE^\gam$ but with
$\psi(s)$ replaced by $-h^{-1}cs^2$. We may write
$\mu=\mu_1+\mu_2$ where
\begin{eqnarray*}
\mu_1(s)&=&h^{-1/4}\{ \psi_{-1}(s)- h^{-1}( i\xi
s+ks^2/2)\}\chi(s)\exp\{-h^{-1}cs^2\}\\
\mu_2(s)&=&h^{-1/4}\left(
\sum_{m=0}^nh^m\psi_{m}(s)\right)\chi(s)\exp\{-h^{-1}cs^2\}
\end{eqnarray*}
Since
\[
| \psi_{-1}(s)- h^{-1}( i\xi s+ks^2/2)|\leq c_{-1}h^{-1}|s|^3
\]
we have $\mu_1\in\tilde{\cE}^1$. Since $|\psi_m(s)|\leq c_m|s|$
for all $s$ we also have $\mu_2\in\tilde{\cE}^1$. The estimate of
this lemma now follows by an obvious modification of
Lemma~\ref{asympt}.


\section{A Semi-classical Transform\label{transform}}

We continue with the assumptions and notation of the last section.
Theorem~\ref{mainasympt} provides the information needed for the
application of Theorem~\ref{mainseq}. We define the set $\Ome$ in
Theorem~\ref{mainseq} by (\ref{nontriv}) and take $\sig$ to be the
semi-classical principal symbol (\ref{scps}) of $A$. In numerical
applications, one would, of course, have to restrict to a finite
subset of $\Ome$, as described in \cite{EBD2}.

We fix $n$ and put $e_{h,u,\xi}=f_{h,u,\xi}/\norm
f_{h,u,\xi}\norm$ where $f_{h,u,\xi}=f$ is defined by
(\ref{fdef}). The semiclassical integral transform $E:L^1(\Ome)\to
L^2(\R)$ is then defined by
\[
(E_h\phi)(x)=\int_\Ome \phi(u,\xi)e_{h,u,\xi}(x)\,\rmd u\rmd \xi
\]
and has norm $1$ by \cite[Theorem VI.8.6]{DS}. The functions
$f_{h,u,\xi}(x) $ are very complicated for large $n$, and the
following approximation may therefore be valuable.

\begin{theorem} Given $h,u,\xi$, let
\begin{equation}
g_{h,u,\xi}(x)=h^{-1/4}\exp\{h^{-1}(i\xi(x-u)+k_{u,\xi}(x-u)^2/2)\}\label{gdef}
\end{equation}
where
\begin{equation}
k_{u,\xi}=-i\frac{\partial \sig}{\partial u}\left\{ \frac{\partial
\sig}{\partial \xi}\right\}^{-1}.\label{kdef}
\end{equation}
If $(u,\xi)\in\Ome$ then $g_{h,u,\xi}\in L^2(\R)$. Define
$E_h^\pr: L^1(\Ome)\to L^2(\R)$ by
\begin{equation}
(E_h^\pr\phi)(x)=\int_\Ome \phi(u,\xi)e^\pr_{h,u,\xi}(x)\,\rmd
u\rmd \xi  \label{Ehprdef}
\end{equation}
where $e^\pr_{h,u,\xi}=g_{h,u,\xi}/\norm g_{h,u,\xi}\norm$. Then
$\norm E_h^\pr \norm=1$ and
\begin{equation}
\lim_{h\to 0}\norm E_h\phi-E_h^\pr\phi\norm =0 \label{strongconv}
\end{equation}
for all $\phi\in L^1(\Ome)$. If we replace $\Ome$ by a compact
subset $U$ of $\Ome$ then
\begin{equation}
\lim_{h\to 0}\norm E_h-E_h^\pr\norm =0 \label{normconv}
\end{equation}
\end{theorem}

\Proof We start by observing that  $\Re(k_{u,\xi})<0$ if and only
if $(u,\xi)\in\Ome$, so $g_{h,u,\xi}\in L^2(\R)$ under the same
conditions. We have $\norm E_h^\pr \norm =1$ by \cite[Theorem
VI.8.6]{DS}.

Let $\{\Ome_n\}_{n=1}^\infty$ be an increasing sequence of compact
subsets of $\Ome$ whose union equals $\Ome$. If we can prove that
the restrictions $E_{h,n}$ and $E_{h,n}^\pr$ to $L^1(\Ome_n)$
satisfy
\begin{equation}
\lim_{h\to 0}\norm E_{h,n}-E_{h,n}^\pr\norm =0\label{normlim}
\end{equation}
then (\ref{strongconv}) and (\ref{normconv}) follow by  standard
procedures.

In Lemma~\ref{gaussapp} we proved that
\[
\norm f_{h,u,\xi}-g_{h,u,\xi}\norm =O(h^{1/2})
\]
for each $(u,\xi)\in\Ome$ as $h\to 0$. The dependence of the error
upon $u,\xi$ and the coefficients of $A$ was given explicitly, and
implies that
\[
\lim_{h\to 0} \sup\{  \norm f_{h,u,\xi}-g_{h,u,\xi}
\norm:(u,\xi)\in\Ome_n\} =0.
\]
Taking (\ref{eq1}) into account we deduce that
\[
\lim_{h\to 0} \sup\{  \norm e_{h,u,\xi}-e^\pr_{h,u,\xi}
\norm:(u,\xi)\in\Ome_n\} =0.
\]
This implies (\ref{normlim}).

\begin{lemma}
Let let $E^\pr_{h,U}$ denote the restriction of $E_h^\pr$ to the
subset $U$ of $\Ome$. If $U,V$ are two compact subsets of $\Ome$
which are spatially disjoint in the sense that $(u,\xi)\in U$ and
$(v,\eta)\in V$ implies $u\not= v$ then the ranges of
$E^\pr_{h,U}$ and $E^\pr_{h,V}$ are uniformly asymptotically
orthogonal in the sense that
\[
\lim_{h\to 0} \norm (E^\pr_{h,U})^\ast E^\pr_{h,V}\norm =0.
\]
The convergence is exponentially fast.
\end{lemma}

\Proof Let $W$ be an open subset of $\R$ such that
$U\subseteq(W\times \R)$ and
$V\cap(\overline{W}\times\R)=\emptyset$. Let $P$ be the orthogonal
projection in $L^2(\R)$ whose range consists of all functions with
support in $W$. Then
\begin{eqnarray*}
\norm (E^\pr_{h,U})^\ast E^\pr_{h,V}\norm&\leq& \norm
(E^\pr_{h,U})^\ast (I-P)E^\pr_{h,V}\norm+\norm (E^\pr_{h,U})^\ast PE^\pr_{h,V}\norm\\
&\leq& \norm (I-P)E^\pr_{h,U}\norm+\norm PE^\pr_{h,V}\norm.
\end{eqnarray*}
We consider further only the first term on the RHS; the other is
treated in a similar manner. If $\phi\in L^1(U)$ then
\begin{eqnarray*}
\norm (I-P)E^\pr_{h}\phi\norm&=&\norm\int_U
(I-P)e^\pr_{h,u,\xi}\phi(u,\xi)\,\rmd u\rmd\xi\norm\\
&\leq& \int_U \norm (I-P)e^\pr_{h,u,\xi}\norm \, |\phi(u,\xi)|\,\rmd u\rmd\xi\\
&\leq& \sup\{\norm (I-P)e^\pr_{h,u,\xi}\norm :(u,\xi)\in
U\}\,\norm\phi\norm\\
 &\leq & \frac{\sup\{\norm
(I-P)g_{h,u,\xi}\norm :(u,\xi)\in U\} }{\inf\{\norm
g_{h,u,\xi}\norm :(u,\xi)\in U\}}\norm\phi\norm.
\end{eqnarray*}
The explicit expression (\ref{gdef}) for $g$ and the compactness
of $U$ ensure that the final supremum converges to $0$
exponentially fast as $h\to 0$ while the final infimum converges
to a positive limit.

If we subdivide $\R$ into small intervals then the lemma implies
that $E^\pr_h$ (or more exactly its restriction to any compact
subregion of $\Ome$) acts asymptotically independently on
subintervals which are not adjacent. If each interval is small
enough we may approximate $E^\pr_h$ in any subinterval by the
operator with a frozen value of $u$.

We conjecture that under suitable conditions on the coefficients
of $A$ the transforms $E_h$ and $E^\pr_h$ are both bounded from
$L^2(\Ome)$ to $L^2(\R)$. As evidence for this we treat the case
in which the variable $u$ in $k_{u,\xi}$ is frozen at the value
$v$. We also assume that $A$ is a \Schrodinger operator, so that
its symbol is of the form $\sig(u,\xi)=\xi^2+c(u)$. This implies
that $k_{v,\xi}=-1/\kap\xi$ where $\kap=2/i c^\pr(v)$. Assuming
that $\kap$ has positive real part, it is immediate that $\Re \,
k_{v,\xi}<0$ if and only if $\xi
>0$. We therefore put $\R^2_{+}=\{ (u,\xi):u\in\R,\xi
>0\}$.

We define the distorted FBI transform $\tilde{E}_{h}:
C_c(\R^2_{+})\to L^2(\R)$ by
\begin{equation}
\tilde{E}_h\phi=h^{-1/2}\int_{\R^2_{+}}
\phi(u,\xi)\tilde{e}_{h,u,\xi}\,\rmd u\rmd \xi \label{Etildedef}
\end{equation}
where $\tilde{e}_{h,u,\xi}=\tilde{g}_{h,u,\xi}/\norm
\tilde{g}_{h,u,\xi}\norm$ and
\begin{equation}
\tilde{g}_{h,u,\xi}(x)=\exp\{i\xi(x-u)/h-(x-u)^2/2h\kap\xi\}.\label{gtildedef}
\end{equation}

\begin{theorem} \label{mainasymptot}
If $\Re( \kap) >0$ and $h>0$ then the operator
(\ref{Etildedef}) may be extended to a bounded operator from
$L^2({\R^2_{+}})$ to $L^2(\R)$ whose norm is bounded above
uniformly as $h\to 0$.
\end{theorem}

\Proof In this proof we write $c_r$ to denote positive constants
which depend only on $\kap$. We always take $\xi$ to be positive.
We have
\begin{eqnarray*}
\norm \tilde{g}_{h,u,\xi}\norm^2&=&  \int_\R \exp\{
-\Re(1/\kap)(x-u)^2/h\xi\}\,\rmd x\\
&=& c_1 h^{1/2}\xi^{1/2}.
\end{eqnarray*}
Therefore
\[
\norm \tilde{g}_{h,u,\xi}\norm =c_2h^{1/4}\xi^{1/4}.
\]

We prove the $L^2$ boundedness of $\tilde{E}^\ast_h$ rather than
that of $\tilde{E}_h$. We have
\begin{equation}
(\tilde{E}^\ast_h f)(u,\xi)= \int_\R K(u,\xi,h,x)f(x)\,\rmd
x\label{kernel}
\end{equation}
where
\begin{eqnarray*}
K(u,\xi,h,x)&=& h^{-1/2}\overline{\tilde{e}_{h,u,\xi}(x)}=\bet_{h,\xi}\gam_\xi(u-x),\\
\bet_{h,\xi}&=& c_3 h^{-3/4}\xi^{-1/4},\\
\gam_\xi(u)&=&\exp\{ i\xi u/h- u^2/2h\overline{\kap}\xi \}.
\end{eqnarray*}
We next take the Fourier transform $\cF$ of (\ref{kernel}) in the
$u$ variable, noting that $\cF$ is a unitary operator on
$L^2({\R^2_{+}})$. This yields
\[
\norm \tilde{E}^\ast_h f \norm =\norm k \norm
\]
where
\[
k(s,\xi)=\bet_{h,\xi}\hat{\gam}_\xi(s)(\cF f)(s)
\]
and
\begin{eqnarray*}
\hat{\gam}_\xi(s)&=&\int_\R \exp\{ iu(\xi/h-s)-u^2/2h\overline{\kap}\xi \}\,\rmd u\\
&=& c_4h^{1/2}\xi^{1/2}\exp\{-(\xi/h-s)^2h\overline{\kap} \xi/2\}.
\end{eqnarray*}
We deduce that
\[
\norm \tilde{E}^\ast_h f \norm\leq c_5\norm \cF f\norm =c_5\norm
f\norm
\]
for all $f\in L^2(\R)$ if and only if
\[
\sup_{s\in\R} \left\{ \int_0^\infty \left|
\bet_{h,\xi}\hat{\gam}_\xi(s)\right|^2 \,\rmd \xi \right\}\leq c_5^2.
\]
Our task therefore, is to prove that the function
\begin{equation}
F(h,s)=\int_0^\infty
h^{-1/2}\xi^{1/2}\exp\{-c_6(\xi/h-s)^2h\xi\}\,\rmd\xi\label{Fhsdef}
\end{equation}
is bounded on $\R^+\times \R$, provided $c_6>0$. If $s\leq 0$ then
putting $\xi=h^{1/3}\eta$ we obtain
\begin{eqnarray*}
F(h,s)&\leq & F(h,0)\\
&=& \int_0^\infty
h^{-1/2}\xi^{1/2}\exp\{-c_6\xi^3/h \}\,\rmd\xi\\
&=&  \int_0^\infty \eta^{1/2}\exp\{-c_6\eta^3 \}\,\rmd\eta,
\end{eqnarray*}
which is finite. If $s>0$ then putting $\xi=hs\eta$ we obtain
\begin{equation}
F(h,s)=G(h^2s^3)\label{FequalsG}
\end{equation}
where
\[
G(t)=\int_0^\infty \eta^{1/2}t^{1/2}\rme^{-c_6(\eta-1)^2\eta t}
\rmd \eta,
\]
so we have to prove that $G$ is bounded on $(0,\infty)$. We do
this in stages. If $0\leq t\leq 1$ then
\[
\int_0^{1/2} \eta^{1/2}t^{1/2}\rme^{-c_6(\eta-1)^2\eta t} \rmd
\eta \leq 1/2
\]
because every term in the integrand is less than $1$. If $t\geq 1$
then
\begin{eqnarray*}
\int_0^{1/2} \eta^{1/2}t^{1/2}\rme^{-c_6(\eta-1)^2\eta t} \rmd
\eta
&\leq& \int_0^{1/2} \eta^{1/2}t^{1/2}\rme^{-c_6\eta t/4} \rmd \eta\\
&\leq& \int_0^\infty \eta^{1/2}t^{1/2}\rme^{-c_6\eta t/4} \rmd \eta\\
&=& c_7 t^{-1}\leq c_7.
\end{eqnarray*}
If $t>0$ then
\begin{eqnarray*}
\int_{1/2}^{4} \eta^{1/2}t^{1/2}\rme^{-c_6(\eta-1)^2\eta t} \rmd
\eta
&\leq& \int_{1/2}^4 2t^{1/2}\rme^{-c_6(\eta-1)^2 t/2} \rmd \eta\\
&\leq& \int_{-\infty}^\infty 2t^{1/2}\rme^{-c_6\eta^2 t/2} \rmd \eta\\
&=& c_8.
\end{eqnarray*}
Finally if $t>0$ then putting $\eta=\zeta t^{-1/3}$ we obtain
\begin{eqnarray*}
\int_{4}^\infty \eta^{1/2}t^{1/2}\rme^{-c_6(\eta-1)^2\eta t} \rmd
\eta
&\leq& \int_{4}^\infty \eta^{1/2}t^{1/2}\rme^{-c_8\eta^3 t} \rmd \eta\\
&\leq& \int_0^\infty \zeta^{1/2}\rme^{-c_8\zeta^3 } \rmd \zeta\\
&=& c_9.
\end{eqnarray*}

One cannot expect $E_h^\ast$ to be isometric, as is the case for
the FBI transform, but we prove that this is asymptotically true
in the semi-classical limit, up to a normalizing constant $c$,
which could be evaluated explicitly.

\begin{theorem} There exists a positive constant $c$ such that
\[
\lim_{h\to 0} \norm E_h^\ast f\norm =c \norm f \norm
\]
for all $f\in L^2(\R)$.
\end{theorem}

\Proof In the proof of Theorem~\ref{mainasymptot} we obtained the
formula
\[
\norm E_h^\ast f\norm^2=c_{10} \int _{-\infty}^\infty F(h,s)|(\cF
f)(s)|^2\,\rmd s
\]
where
\[
0\leq F(h,s)\leq c_{11}
\]
for all $h>0$ and $s\in\R$. By the dominated convergence theorem
it suffices to prove that
\[
\lim_{h\to 0} F(h,s)=c_{12}:=\int_0^\infty
\eta^{1/2}\exp\{-c_6\eta^3\}\,\rmd \eta
\]
for all $s\in \R$. We do this for $s>0$, noting that the cases
$s=0$ and $s<0$ are similar. By (\ref{FequalsG}) it suffices to
prove that $\lim_{t\to 0+}G(t)=c_{12}$. As $t\to 0+$ we have
\begin{eqnarray*}
G(t)&= & \int_0^\infty \eta^{1/2}t^{1/2}\rme^{-c_6(\eta-1)^2\eta
t} \rmd \eta\\
&\sim & \int_0^\infty \eta^{1/2}t^{1/2}\rme^{-c_6\eta^3 t} \rmd
\eta\\
&=& c_{12}
\end{eqnarray*}
using the change of variable $\eta\to \eta t^{-1/3}$.

In order to extend Theorem~\ref{mainasymptot} to second order
differential operators other than \Schrodinger operators, it needs
to be generalized as follows.

\begin{theorem}
Let $\kap:(0,\infty)\to\C$ be a continuous function, let
$c_0,c_\infty$ be positive constants and let $\alp_0,\alp_\infty$
be non-negative constants such that
\[
\begin{array}{rccclcl}
c_0^{-1}\xi^{\alp_0} &\leq &\Re\kap(\xi)&\leq & c_0\xi^{\alp_0}
&\mbox{ if}& 0<\xi\leq 1,\\
c_\infty^{-1}\xi^{\alp_\infty} &\leq &\Re\kap(\xi)&\leq &
c_\infty\xi^{\alp_\infty} &\mbox{ if}& 1\leq \xi <\infty.\\
\end{array}
\]
Then the conclusion of Theorem~\ref{mainasymptot} is still valid
if we replace (\ref{gtildedef}) by
\[
\tilde{g}_{h,u,\xi}(x)=\exp\{i\xi(x-u)/h-(x-u)^2/2h\kap(\xi) \}.
\]
\end{theorem}

\Proof We make obvious adaptations to the proof of
Theorem~\ref{mainasymptot} up to (\ref{Fhsdef}), which becomes
\begin{eqnarray*}
F(h,s) &=& \int_0^\infty
h^{-1/2}(\Re\kap(\xi))^{1/2}\exp\{-c_6(\xi/h-s)^2h\,\Re\kap(\xi)\}\,\rmd\xi \\
&\leq & \int_0^1
h^{-1/2}c_0^{1/2}\xi^{\alp_0/2}
\exp\{-c_6(\xi/h-s)^2hc_0^{-1}\xi^{\alp_0}\}\,\rmd\xi \\
&& +\int_1^\infty
h^{-1/2}c_\infty^{1/2}\xi^{\alp_\infty/2}
\exp\{-c_6(\xi/h-s)^2hc_\infty^{-1}\xi^{\alp_\infty}\}\,\rmd\xi.
\end{eqnarray*}
Each of these integrals is estimated by the same method as in
Theorem~\ref{mainasymptot}.


\section{Constructing the boundary pseudospectra\label{bdrsect}}

When one examines the pseudo-eigenfunctions in several exactly
soluble examples, \cite{EBD3,EBD2,red,RT}, one sees that they do
not conform to the above ideas. They are strongly localized at one
end of the interval in question, and decrease exponentially as one
moves away from this end.

In this section we develop the general theory of boundary
pseudospectra for variable coefficient operators in the
one-dimensional context. A partial extension to higher dimensions
and manifolds is described in the next section. We assume that
\[
(L_hf)(x)=-h^2a(x)f^{\pr\pr}(x)-ihb(x)f^\pr(x)+c(x)f(x)
\]
for $x\in [0,\gam]$.  The semiclassical principal symbol is
\[
\sig(u,\xi)=a(u)\xi^2+b(u)\xi +c(u).
\]
We will need the fact that the symbol can be analytically
continued to complex $\xi$, but only assume the coefficients of
$L_h$, and therefore $\sig$, to be $C^\infty$ in $u$ on
$[0,\gam]$. Similar but weaker estimates can be proved if the
coefficients are only $C^n$ for some $n$. We assume ellipticity,
in other words that $a(x)\not= 0$ for all $x\in[0,\gam]$. We start
by ignoring the boundary conditions and looking for a
pseudo-eigenfunction of the form
\begin{equation}
f(s)=h^{-1/2}\chi(s)\exp(\psi(s))\label{bdryf}
\end{equation}
where
\[
\psi(s)=\sum_{m=-1}^n h^m\psi_m(s).
\]

We assume that $\chi\in C^\infty[0,\gam]$ satisfies $\chi(s)=1$ if
$0\leq s\leq \del/2$ and $\chi(s)=0$ if $s\geq \del$; the constant
$\del
>0$ must be small enough for the proof of Theorem~\ref{main2} to
be valid. We put
\[
\psi_{-1}(s)=i\int_{v=0}^s\left\{
-\frac{b(v)}{2a(v)}+\sqrt{w(\xi,v)}\right\}\,\rmd v
\]
where
\[
w(\xi,v)= \frac{a(0)\xi^2}{a(v)} + \frac{b(0)\xi}{a(v)}
+\frac{b(v)^2}{4a(v)^2} +\frac{c(0)-c(v)}{a(v)}.
\]
As before we take the branch of the square root which equals
$\xi+b(0)/2a(0)$ at $v=0$. However we now require $\Im(\xi)>0$, in
order to ensure that $f(s)$ decays rapidly as $s$ increases. We
have
\[
\psi_{-1}(s)=i\xi s+ks^2/2+O(s^3)
\]
for small $s>0$ as before.

\begin{lemma} \label{asympt2}
Let $F$ be a positive continuous function on $[0,\del]$ and let
$G$ be a continuous function on $[0,\del]$. If $m$ is a
non-negative even integer then
\[
\int_{0}^\del s^mG(s)\exp\{-h^{-1}sF(s)\}\,\rmd s\sim ch^{m+1}
\]
as $h\to 0+$, where
\[
c= \frac{G(0)\Gam(m+1)}{ F(0)^{m+1}}.
\]
\end{lemma}

In the following theorem we put $(Qf)(x)=xf(x)$ and
$(Pf)(x)=-ihf^\pr(x)$ as before. Although $Q$ is self-adjoint on
an obvious domain, we impose no boundary conditions on $P$, which
is therefore not even symmetric.

\begin{theorem}\label{main2}
If the coefficients of $L_h$ are $C^\infty$ and $\Im(\xi)>0$ then
for any positive integer $n$ there exist functions $f$ which
depend on $h,n,\xi$ such that
\begin{eqnarray}
\lim_{h\to 0}\norm {f}\norm&= & c>0\label{eq1b}\\
\norm Q{f}\norm&= & O(h)\label{eq2b}\\
\norm P{f}-\xi{f}\norm&= & O(h)\label{eq3b}\\
\norm L_h{f}-\sig(0,\xi){f}\norm &=& O(h^{n+2})\label{eq4b}
\end{eqnarray}
as $h\to 0$.
\end{theorem}

\Proof Let $f$ be given by (\ref{bdryf}). To prove (\ref{eq1b}) we
write
\begin{eqnarray*}
\norm {f}\norm^2 & = & h^{-1}\int_{0}^\del
\chi(s)^2\exp\{  2\Re(\psi(s))\}\,\rmd s\\
&=& h^{-1}\int_{0}^\del G(s)\exp\{-h^{-1}sF(s)\}\,\rmd s
\end{eqnarray*}
where
\begin{eqnarray*}
F(s)& = & - 2\Re(\psi_{-1}(s))/s \\
G(s) &=&\chi(s)^2\exp\left\{ 2\Re\left( \sum_{m=0}^n h^m
\psi_m(s). \right) \right\}
\end{eqnarray*}
This is of the form treated by Lemma~\ref{asympt2} if $\del>0$ is
small enough to ensure that $F(s)>0$ for all $s\in [0,\del]$.

To prove (\ref{eq2b}) we write
\[
\norm {Qf}\norm^2 = h^{-1}\int_{0}^\del
s^2G(s)\exp\{-h^{-1}sF(s)\}\,\rmd s
\]
and apply Lemma~\ref{asympt2} again.

The proof of (\ref{eq3b}) uses Lemma~\ref{asympt2} and the
expansion
\[
Pf-\xi f=\mu_1+\mu_2+\mu_3
\]
where
\begin{eqnarray*}
\mu_1 &=& -ih^{-1/2}\{ \psi_{-1}^\pr(s)- i\xi)\}\chi(s)\exp\{\psi(s)\}\\
\mu_2 &=& -ih^{1/2}\left(
\sum_{m=0}^nh^m\psi_{m}^\pr(s)\right) \chi(s)\exp\{\psi(s)\}\\
\mu_3 &=& -ih^{1/2} \chi^\pr(s)\exp\{\psi(s)\}.
\end{eqnarray*}

The proof of (\ref{eq4b}) follows in a similar way from the
formula
\[
L_h{f}-\sig(0,\xi){f}=\left( \sum_{m=n+2}^{2n+2} h^m \phi_m\right)
f+O(h^\infty).
\]

We finally assume the boundary conditions
\begin{equation}
uh f^\pr(0)+w f(0)=0\label{bdry}
\end{equation}
for some complex constants $u,\, w$, not both zero. We say that
$L_h$ satisfies the exit condition at $0$ if $\Im(-b(0)/a(0))>0$.
This language is motivated by the example discussed in
\cite{EBD2}, in which $L_h$ is the generator of a subMarkov
diffusion on an interval. Given the exit condition at $0$, we
define the boundary semiclassical pseudospectrum at $0$ to be the
set
\begin{equation}
\tilde{\Lam}=\{ \xi:0<\Im(\xi)<\Im(-b(0)/a(0)) \}.\label{neg}
\end{equation}
If $\xi_1\in\tilde{\Lam}$ and $z=\sig(0,\xi_1)$ then the other
solution $\xi_2$ of $\sig(0,\xi)=z$ also lies in $\tilde{\Lam}$.
We have $\xi_1=\xi_2$ if and only if $z= c(0)-b(0)^2/4a(0)$. The
set $\sig(0,\tilde{\Lam})$ is the region inside the parabola $P=\{
\sig (0,t):t\in\R\}$.

Those familiar with \cite{EBD3,red,RT} will observe the close
relationship between the above and the winding number calculations
there. At a qualitative level the given operator can be
approximated near the end of the interval by the operator whose
coefficients are frozen to the values which they have at the
endpoint. Our theorem below provides quantitative flesh to this
idea. It also provides the precise form of the relevant
pseudo-eigenfunction, which is not easy to guess from the constant
coefficient case.

\begin{theorem}\label{eight}
Let $L_h$ satisfy the exit condition at $0$ and let $z$ lie inside
the parabola $P$. Assuming $z\not= c(0)-b(0)^2/4a(0)$, let
$\xi_1,\,\xi_2\in\tilde{\Lam}$ denote the two distinct solutions
of $\sig (0,\xi)=z$. Given $h>0$ and $n \geq 1$, let $f_r$ be the
boundary pseudo-eigenfunctions associated with $h,n,\xi_r$ as in
(\ref{bdryf}) and Theorem~\ref{main2}, and let
\begin{equation}
f=(iu\xi_2+w)f_1-(iu\xi_1+w)f_2.\label{funct}
\end{equation}
Then $f$ satisfies the boundary condition (\ref{bdry}) at $0$ and
\begin{equation}
\norm L_h f-zf\norm / \norm f\norm =O(h^{n+2})\label{estt}
\end{equation}
as $h\to 0$.
\end{theorem}

\Proof The assumptions imply that $f_r$ satisfy the estimates of
Theorem~\ref{main2}, from which (\ref{estt}) follows. The proof
that $f$ satisfies (\ref{bdry}) depends upon the identities
$f_r(0)=h^{-1/2}$ and $f_r^\pr(0)=ih^{-3/2}\xi_r$.

\section{Higher Dimensions}

The extension of the above ideas to higher dimensions needs more
machinery. We are mainly interested in bounded regions in $\R^N$
with smooth boundary, but since the proof of our main result
depends upon choosing local coordinates around a boundary point
rather carefully, we write down the argument in a manifold
context. Let $X$ be a smooth $N$-dimensional manifold with
boundary $\partial X$. Let $X$ be provided with a volume measure
$\rmd vol$ which has positive $C^\infty$ density $v(x)$ when
restricted to any coordinate neighbourhood $U$.

The natural differential $\rmd :C^n(X)\to C^{n-1}(T^\ast X)$ is
given within $U$ by
\[
\rmd f(x)=(\partial_1f(x),..., \partial_nf(x))
\]
and the adjoint operator $\rmd^\ast:C^n(TX)\to C^{n-1}(X)$  acts
on a section $g\in C^n(TU)$ by
\[
\rmd^\ast g(x)=-v(x)^{-1}\partial_j(v(x)g^j(x)).
\]

The differential operator $L_h$ is determined by three coefficient
functions, all assumed to be $C^\infty$ and complex-valued on $X$;
we write $T_x$ and $T^\ast_x$ in place of $T_x\otimes \C$ and
$T^\ast_x\otimes \C$ below. We assume that $a(x):T^\ast_x\to T_x$,
$b(x)\in T_x$ and $c(x)\in\C$ for all $x\in X$. Given $h>0$ and
$f\in C^\infty(X)$ we then put
\[
(L_hf)(x)=h^2\rmd^\ast(a(x)\rmd f(x))-ihb(x)\!\cdot\! \rmd
f(x)+c(x)f(x).
\]
Throughout this section a dot indicates the natural action of a
covector on a tangent vector at some point of $X$. In the
coordinate neighbourhood $U$ the above formula may be written in
the form
\[
(L_hf)(x)=-h^2v^{-1}(x)\partial_j\left(v(x)a^{j,k}(x)\partial_kf(x)\right)-ihb^j(x)
\partial_jf(x)+c(x)f(x)
\]
using the usual summation convention, or in the form
\begin{equation}
(L_hf)(x)=-h^2a^{j,k}(x)\partial_{j,k} f(x)-ihb^j(h,x)\partial_j
f(x)+c(h,x)f(x)\label{leftform}
\end{equation}
where
\begin{eqnarray}
b^j(h,x)&=& b^j(x)+hb^j_1(x)\label{bform}\\
c(h,x)&=&c(x)+hc_1(x)+h^2c_2(x)\label{cform}.
\end{eqnarray}
The set of all operators of the form (\ref{leftform}) is invariant
under changes of local coordinates.

The symbol of $L_h$ is given by
\[
\sig_h(x,\xi)=h^2a^{j,k}(x)\xi_j\xi_k +hb^j(h,x)\xi_j+c(h,x)
\]
which is not an invariant expression: $\rmd^\ast$ and $L_h$ both
depend upon the choice of the density $v$. However the
semiclassical principal symbol
\begin{eqnarray*}
\sig(x,\xi)&=&\lim_{h\to
0}\sig_h(x,h^{-1}\xi)\\
&=&a^{j,k}(x)\xi_j\xi_k
+b^j(x)\xi_j+c(x)\\
&=& a(x)\xi\cdot\xi+b(x)\cdot\xi+c(x)
\end{eqnarray*}
is invariant under changes of local coordinates.

The following theorem is a multi-dimensional `boundary' analogue
of Theorem~\ref{two}. We expect that there is also a
multi-dimensional analogue of Theorem~\ref{eight}. We choose a
point in $\partial X$, label it $p$, and choose a complex
cotangent vector $\xi$ at $p$. We require that $\Im(\xi)$ has zero
dot product with any vector at $p$ which is tangent to $\partial
X$ and positive dot product with any inward pointing vector at
$p$. If $U$ is a coordinate neighbourhood around $p$ we always
assume that $p$ is represented by the point $0\in\R^N$.

\begin{theorem} \label{xyz} Let $L_h$ be of the form (\ref{leftform}) where
all of the coefficients in (\ref{leftform}), (\ref{bform}),
(\ref{cform}) are $C^\infty$ functions on $U$. Let the complex
cotangent vector $\xi$ at $0\in\partial X$ satisfy the conditions
of the last paragraph. Then for every sufficiently small $h>0$
there exists $f_h\in C^\infty(X)$ which vanishes outside a
neighbourhood of $0$ whose radius is of order $h^{1/2}$, and
satisfies
\begin{eqnarray}
\lim_{h\to 0}\norm f_h\norm_2&=& c>0\label{mul1}\\
\norm L_h f_h -\sig(0,\xi)f_h\norm_2 &=&O(h^{1/2})\label{mul2}
\end{eqnarray}
as $h\to 0$.
\end{theorem}

\Proof Let $\R^N_+$ denote the set of $x\in\R^N$ for which
$x^N\geq 0$ and let $\R^N_0$ denote the set of $x$ for which
$x^N=0$. We choose local coordinates around $0$ such that
\[
U=\{ x\in\R^N_+:|x|<\rho\}
\]
and put
\[
\partial U = \{ x\in\R^N_0:|x|<\rho\}
\]
for some $\rho >0$. We write $x=(x^\pr,x^N)$ where
$x^\pr\in\R^{N-1}$ and $x^N\in\R$. Our assumptions imply that
$\xi=(\xi^\pr,\xi^N)$ where $\xi^\pr$ is real and
$\eta:=\Im(\xi^N)>0$.

Put $\alp=1/2$ and $\gam=(N+1)/4$. Let $\phi_1$ be a smooth
function on $\R^{N-1}$ which equals $1$ if $|x^\pr|\leq 1$ and $0$
if $|x^\pr|\geq 2$. Let $\phi_2$ be a smooth function on
$[0,\infty)$ which equals $1$ if $0\leq x^N\leq 1$ and $0$ if
$x^N\geq 2$. Let $\phi(x)=\phi_1(x^\pr)\phi_2(x^N)$. Then the
smooth function
\[
f_h(x)=h^{-\gam}\rme^{ih^{-1}\xi\cdot x}\phi(h^{-\alp}x)
\]
on $U$ has support with the required property for all small enough
$h>0$.

To prove (\ref{mul1}) we observe that
\begin{eqnarray*}
\norm f_h\norm_2^2&\sim & v(0)h^{-2\gam}
\int_{\R^{N-1}}\phi_1(h^{-\alp}x^\pr)^2\,\rmd^{N-1}
x^\pr\int_0^\infty \rme^{-2h^{-1}\eta x^N} \phi_2(h^{-\alp}
x^N)^2\,\rmd x^N\\
&=&
 v(0)h^{-2\gam+(N-1)\alp+1}\int_{\R^{N-1}}\phi_1(y^\pr)^2\,\rmd^{N-1}
y^\pr\int_0^\infty \rme^{-2\eta s} \phi_2(h^{1-\alp}s)^2\,\rmd s\\
&\to& v(0)(2\eta)^{-1}
\int_{\R^{N-1}}\phi_1(y^\pr)^2\,\rmd^{N-1}y^\pr
>0
\end{eqnarray*}
as $h\to 0$.

The proof of (\ref{mul2}) depends upon writing
\[
L_h f_h-\sig(0,\xi)f_h=\sum_{m=1}^7 g_m
\]
where
\begin{eqnarray*}
g_1&=& h^{-\gam}\{a^{j,k}(x)-a^{j,k}(0)\}
\xi_j\xi_k\rme^{ih^{-1}\xi\cdot x}\phi(h^{-\alp}x)\\
g_2&=&-ih^{1-\alp-\gam}a^{j,k}(x)\xi_j\rme^{ih^{-1}\xi\cdot
x}\phi_k(h^{-\alp}x)\\
g_3&=&-ih^{1-\alp-\gam}a^{j,k}(x)\xi_k\rme^{ih^{-1}\xi\cdot
x}\phi_j(h^{-\alp}x)\\
g_4&=&-h^{2-2\alp-\gam}a^{j,k}(x)\rme^{ih^{-1}\xi\cdot
x}\phi_{j,k}(h^{-\alp}x)\\
g_5&=& h^{-\gam}\{ b^j(h,x)-b^j(0)\}\xi_j\rme^{ih^{-1}\xi\cdot x}\phi(h^{-\alp}x)\\
g_6&=& -ih^{1-\alp-\gam}b^j(h,x)\rme^{ih^{-1}\xi\cdot
x}\phi_j(h^{-\alp}x)\\
g_7&=& h^{-\gam} \{ c(h,x)-c(0)\} \rme^{ih^{-1}\xi\cdot
x}\phi(h^{-\alp}x).\\
\end{eqnarray*}
We estimate the $L^2$ norm of each of these as above, obtaining
$\norm g_r \norm_2=O(h^{\alp})$ for $r=1,5,7$, $\norm g_r
\norm_2=O(h^{1-\alp})$ for $r=2,3,6$ and $\norm g_r
\norm_2=O(h^{2-2\alp})$ for $r=4$. Given these estimates, the
optimal value of $\alp$ is $1/2$.

We next impose boundary conditions of the form
\[
h\,u(x^\pr)n(x^\pr,0)\cdot \rmd f(x^\pr,0)+w(x^\pr)f(x^\pr,0)=0
\]
for all $x^\pr\in\partial U$, where the complex-valued
coefficients $u,\, w$ are $C^\infty$ on $\partial U$; we assume
non-degeneracy of the boundary conditions at $0$ in the sense that
$u(0^\pr)$ and $w(0^\pr)$ do not both vanish. The vector field $n$
on $U$ is supposed to be smooth and transversal in the sense that
it has a non-zero inward pointing component at every point of
$\partial U$. We use the associated flow to construct local
coordinates. In other words we choose local coordinates for which
the boundary conditions can be written in the form
\begin{equation}
h\,u(x^\pr)\partial_N
f(x^\pr,0)+w(x^\pr)f(x^\pr,0)=0.\label{bdrycond}
\end{equation}

We say that the complex covector $\xi$ at $0$ is admissible under
the following conditions. We require that $\Im(\xi)$ has positive
dot product with any inward pointing vector at $0$. We require
that the same conditions hold for a second complex covector
$\tilde{\xi}$ at $0$. We require that
$z:=\sig(0,\xi)=\sig(0,\tilde{\xi})$ and that $\xi\cdot
t=\tilde{\xi}\cdot t\in\R$ for any vector $t$ which is tangent to
$\partial U$ at $0$. In the local coordinates specified above we
are fixing $\xi^\pr=\tilde{\xi}^\pr\in \R^{N-1}$ and assuming that
the two solutions $\xi_N$ and $\tilde{\xi}_N$ of the quadratic
equation
\[
\sig(0,(\xi^\pr,s))=z
\]
in $s\in\C$ both have positive imaginary parts. We say that $L_h$
satisfies the exit condition at $0$ if the set of admissible $\xi$
is non-empty.

\begin{theorem} If $\xi\in\C^N$ is an admissible covector and
$z=\sig(0,\xi)$ then under the above conditions there exist
$g_h\in C^\infty(U)$ satisfying the boundary conditions
(\ref{bdrycond}) and also
\begin{eqnarray}
\supp(g_h)&\subseteq& \{ x\in U:|x|<c^\pr h^{1/2}\}\label{multi1}\\
\lim_{h\to 0}\norm g_h\norm_2&=& c>0\label{multi2}\\
\norm L_h g_h -zg_h\norm_2 &=&O(h^{1/2})\label{multi3}
\end{eqnarray}
as $h\to 0$.
\end{theorem}

\Proof We put
\[
g_h(x)=\alp(x^\pr)f_h(x)+\tilde{\alp}(x^\pr)\tilde{f}_h(x)
\]
where
\begin{eqnarray*}
f_h(x)&=& h^{-\gam}\rme^{ih^{-1}\xi\cdot x}\phi(h^{-1/2}x) \\
\tilde{f}_h(x)&=& h^{-\gam}\rme^{ih^{-1}\tilde{\xi}\cdot
x}\phi(h^{-1/2}x).
\end{eqnarray*}
In this equation $\gam=(N+1)/4$ and $\phi$ is as in the proof of
Theorem~\ref{xyz}. Also $\xi=(\xi^\pr,\xi_N)$ and
$\tilde{\xi}=(\xi^\pr,\tilde{\xi}_N)$. The coefficients $\alp,\,
\tilde{\alp}$ are to be determined. Before continuing, we mention
that in the case of Dirichlet boundary conditions we put
$g_h=f_h-\tilde{f}_h$, that is
$\alp(x^\pr)=-\tilde{\alp}(x^\pr)=1$; most of the calculations
below are much simpler in this situation.

It is immediate from the definition that
\begin{eqnarray*}
g_h(x^\pr,0)&=& h^{-\gam}\{
\alp(x^\pr)+\tilde{\alp}(x^\pr)\}\rme^{ih^{-1}\xi^\pr\cdot
x^\pr}\phi_1(h^{-1/2}x^\pr)\\
h\partial_N g_h(x^\pr,0)&=& h^{-\gam}\{
\alp(x^\pr)\xi_N+\tilde{\alp}(x^\pr)\tilde{\xi}_N\}\rme^{ih^{-1}\xi^\pr\cdot
x^\pr}\phi_1(h^{-1/2}x^\pr).
\end{eqnarray*}
It follows that $g_h$ satisfies the boundary conditions provided
\[
iu(x^\pr)\{
\alp(x^\pr)\xi_N+\tilde{\alp}(x^\pr)\tilde{\xi}_N\}+w(x^\pr)\{
\alp(x^\pr)+\tilde{\alp}(x^\pr)\}=0.
\]
This is solved by putting
\begin{eqnarray*}
\alp(x^\pr)&=& w(x^\pr)+iu(x^\pr)\tilde{\xi}_N   \\
\tilde{\alp}(x^\pr)&=&  -w(x^\pr)-iu(x^\pr)\xi_N   .
\end{eqnarray*}
Since $\xi_N\not= \tilde{\xi}_N$, $\alp$ and $\tilde{\alp}$ cannot
both vanish near $0^\pr$.

The validity of (\ref{multi1}) is immediate. To prove
(\ref{multi2}) we note that
\begin{eqnarray*}
\norm g_h\norm_2^2&=& h^{-2\gam}\int_{\R^N} \left|
\alp(x^\pr)\rme^{ih^{-1}\xi\cdot
x}\phi(h^{-1/2}x)+\tilde{\alp}(x^\pr)\rme^{ih^{-1}\tilde{\xi}\cdot
x}\phi(h^{-1/2}x)\right|^2 v(x)\,\rmd^N x\\
&\sim& h^{-2\gam}\int_{\R^N} \left|
\alp(0^\pr)\rme^{ih^{-1}\xi\cdot
x}\phi(h^{-1/2}x)+\tilde{\alp}(0^\pr)\rme^{ih^{-1}\tilde{\xi}\cdot
x}\phi(h^{-1/2}x)\right|^2 v(0)\,\rmd^N x\\
&=& h^{-2\gam}v(0)\int_{\R^{N-1}}
\phi_1(h^{-1/2}x^\pr)^2\,\rmd^{N-1}
x^\pr\times\\
&&\int_0^\infty \left| \alp(0^\pr)\rme^{ih^{-1}\xi_N
x^N}+\tilde{\alp}(0^\pr)\rme^{ih^{-1}\tilde{\xi}_N
x^N}\right|^2 \phi_2(h^{-1/2}x^N)^2\,\rmd x^N\\
&=& v(0)\int_{\R^{N-1}} \phi_1(s^\pr)^2\,\rmd^{N-1}
s^\pr\times\\
&&\int_0^\infty \left| \alp(0^\pr)\rme^{i\xi_N
s^N}+\tilde{\alp}(0^\pr)\rme^{i\tilde{\xi}_N
s^N}\right|^2 \phi_2(h^{1/2}s^N)^2\,\rmd s^N\\
&\to & v(0)\int_{\R^{N-1}} \phi_1(s^\pr)^2\,\rmd^{N-1}
s^\pr\,\int_0^\infty \left| \alp(0^\pr)\rme^{i\xi_N
s^N}+\tilde{\alp}(0^\pr)\rme^{i\tilde{\xi}_N s^N}\right|^2
\,\rmd s^N\\
& >& 0
\end{eqnarray*}
as $h\to 0$.

The proof of (\ref{multi3}) depends upon writing
\[
L_hg_h-zg_h=k_1+k_2+k_3+k_4
\]
where
\begin{eqnarray*}
k_1(x)&=& \alp(0^\pr)\{L_hf_h(x)-zf_h(x)\}        \\
k_2(x)&=&   \tilde{\alp}(0^\pr)\{L_h\tilde{f}_h(x)-z\tilde{f}_h(x)\}   \\
k_3(x)&=&    L_h[\{\alp(x^\pr)-\alp(0^\pr)\}f_h(x)]           \\
k_4(x)&=&L_h[\{\tilde{\alp}(x^\pr)-\tilde{\alp}(0^\pr)\}\tilde{f}_h(x)]
\end{eqnarray*}
and then estimating each term as before.

{\bf Acknowledgements} I should like to thank L N Trefethen and Y
Safarov for helpful comments. I also acknowledge financial support
under the EPSRC grant GR/R81756/01.

\vspace{2 cm}

Department of Mathematics \\
King's College\\
Strand\\
London\\
WC2R 2LS\\
England

E.Brian.Davies@kcl.ac.uk


\vspace{2 cm}

\begin{thebibliography}{99}

\bibitem{boe}  B\"ottcher A: Pseudospectra and singular values of large convolution operators
J. Integral Equations Appl. 6 (1994), 267-301.

\bibitem{QTOS} Davies E B: Quantum Theory of Open Systems. Acad.
Press, London, 1976.

\bibitem{EBD5} Davies E B: Pseudospectra, the harmonic operator and
complex resonances. Proc. Roy. Soc. London A 455 (1999) 585-599.

\bibitem{EBD6} Davies E B: Semi-classical states for
non-self-adjoint \Schrodinger operators. Commun. Math. Phys. 200
(1999) 35-41.

\bibitem{EBD3} Davies E B: Pseudospectra of differential
operators. J. Operator Theory 43 (2000) 243-262.

\bibitem{EBD2} Davies E B: Computing the decay of a simple
reversible sub-Markov semigroup. Preprint, 2003. LMS J Numer.
Anal. to appear.

\bibitem{EBDwebsite} Davies E B: Spectral Theory, chs. 1-3, 2003.
{\tt http://www.mth.kcl.ac.uk/
\\ MAO/Semigroups/DaviesEBrianSemigroupsBook.html}

\bibitem{dsz} Dencker N, Sj\"ostrand J, Zworski M: Pseudospectra
of semi-classical (pseudo)differential operators. Comm. Pure Appl.
Math. to appear.

\bibitem{DT} Driscoll T A, Trefethen L N: Pseudospectra for the wave
equation with an absorbing boundary. J. Comput. Appl. Math. 69
(1996), 125-142.

\bibitem{DS} Dunford N, Schwartz J T: Linear Operators, Part 1.
Interscience Inc., New York, 1966.

\bibitem{ETwebsite}M. Embree and L. N. Trefethen. Pseudospectra
Gateway.\\
{\tt http://www.comlab.ox.ac.uk/pseudospectra}

\bibitem{han} Hansen PC: Rank-Deficient and Discrete Ill-Posed Problems:
Numerical Aspects of Linear Inversion. SIAM Monographs on
Mathematical Modeling and Computation 4, 1997.

\bibitem{HS} Hayashi M, Sakaguchi F: Subnormal operators regarded
as generalized oobservables and compound-system-type normal
extension related to $su(1,1)$. J. Phys. A, Math. Gen., 33 (2000)
7793-7820.

\bibitem{HT} Higham N J, Tisseur F: More on pseudospectra for
polynomial eigenvalue problems and applications in control theory.
Linear Alg. Applic. 351-352 (2002) 435-453.

\bibitem{hol} Holevo A S: Probabilistic and Statistical Aspects
of Quantum Theory. North-Holland, Amsterdam, 1982.

\bibitem{lav} Lavall\'ee P-F: Ph D thesis, Rennes, France, 1997.

\bibitem{mar}Martinez A: An Introduction to Semiclassical and
Microlocal Analysis. Springer, New York, 2002.

\bibitem{red} Reddy S C: Pseudospectra of Wiener-Hopf operators
and constant coefficient differential operators. J. Integral Eqns.
Applic. 5 (1993) 369-4-3.

\bibitem{RT} Reddy S C, Trefethen L N: Pseudospectra of the
convection-diffusion operator. SIAM J. Appl. Math. 54 (6) (1994)
1634-1649.

\bibitem{saf}Safarov Yu: Pseudodifferential operators and linear
connections. Proc. London Math. Soc. (3) 74 (1997) 379-416.

\bibitem{TH} Tisseur F, Higham N J: Structured pseudospectra for
polynomial eigenvalue problems, with applications. SIAM J. Matrix
Anal. Appl. 23 (1) (2001) 187-208.

\bibitem{tre1}Trefethen L N: Pseudospectra of linear
operators. SIAM Review 39 (1997) 383-406.

\bibitem{tre2}Trefethen L N: Computation of pseudospectra. Acta
Numerica. 8 (1999) 247-295.

\bibitem{CT} Trefethen L N, Chapman S J:
Wave packet pseudomodes of twisted Toeplitz matrices Oxford
University Computing Laboratory Numerical Analysis Technical
Report 02/22, December 2002.

\bibitem{zworski} Zworski M: A remark on a paper by E B Davies. Proc.
Amer. Math. Soc. 129 (2001) 2955-2957.

\end{thebibliography}
\end{document}